\documentclass{article}

\usepackage[top=20truemm,bottom=20truemm,left=20truemm,right=20truemm]{geometry}
\usepackage{authblk}
\usepackage{abstract}
\usepackage{amsmath}
\usepackage{amssymb}
\usepackage{hyperref}
\usepackage[dvipdfmx]{graphicx}
\usepackage{color}

\author[1]{Takuya Tsuchiya\thanks{\href{mailto:t-tsuchiya@hi-tech.ac.jp}{\nolinkurl{t-tsuchiya@hi-tech.ac.jp}}}}
\author[2]{Makoto Nakamura}
\affil[1]{%
  Center for Liberal Arts and Sciences,
  Hachinohe Institute of Technology,
  Japan
}
\affil[2]{%
Department of Pure and Applied Mathematics,
Graduate School of Information Science and Technology,
Osaka University,
1-5 Yamadaoka, Suita, Osaka 565-0871, Japan
}

\title{Numerical simulations of semilinear Klein--Gordon
  equation in the de Sitter spacetime
  with structure-preserving scheme}

\begin{document}

\maketitle

\begin{abstract}
  We perform some simulations of the semilinear Klein--Gordon equation in the de
  Sitter spacetime.
  We reported the accurate numerical results of the equation with the
  structure-preserving scheme (SPS) in an earlier publication (Tsuchiya and
  Nakamura in J. Comput. Appl. Math. \textbf{361}: 396--412, 2019).
  To investigate the factors for the stability and accuracy of the numerical
  results with SPS, we perform some simulations with three discretized
  formulations.
  The first formulation is the discretized equations with SPS, the second one is
  with SPS that replaces the second-order difference as the standard
  second-order central difference, and the third one is with SPS that replaces
  the discretized nonlinear term as the standard discretized expression.
  As a result, the above two replacements in SPS are found to be effective for
  accurate simulations.
  On the other hand, the ingenuity of replacing the second-order difference in
  the first formulation is not effective for maintaining the stability of the
  simulations.
\end{abstract}

\section{%
  Introduction
}%
\label{sec:Introduction}

Stable and accurate numerical simulations are necessary for understanding
natural and social phenomena in detail.
To realize this, numerical methods such as discretizations should be in a
mathematically guaranteed format because the numerical errors mainly
occur during the processes of discretizations.
For numerical schemes of partial differential equations, there are several
well-known methods such as the Crank--Nicolson and Runge--Kutta schemes.
However, it is difficult to perform stable and accurate numerical simulations
for nonlinear partial differential equations since there are large numerical
errors and vibrations in the solutions caused by nonlinearlity.
Thus, suitable schemes have been suggested to perform successful simulations.
One of the schemes is the structure-preserving scheme (SPS)
\cite{Furihata-1999-JComputPhys, Furihata-Matsuo-2010}.
This scheme conserves some structures at the continuous level, and thus enables
stable and accurate numerical simulations.

In this paper, we review the discretized equations of the semilinear
Klein--Gordon equation in the de Sitter spacetime with SPS and perform some
simulations to investigate their stability and accuracy.
For investigating the  semilinear Klein--Gordon equation in the de Sitter
spacetime, analytical \cite{Yagdjian-Galstian-2009-CMP, Yagdjian-2009-DCDS,
  Yagdjian-2013-Springer, Nakamura-2014-JMAA, Nakamura-2021-JMP} and numerical
\cite{Yazici-Sengul-OpenPhys-2016, Tsuchiya-Nakamura-2019-JCAM}
research studies have been conducted.
In \cite{Tsuchiya-Nakamura-2019-JCAM}, we reported some accurate numerical
results of the semilinear Klein--Gordon equation with SPS.
There are some differences between the standard discretized equation and the
discretized equation with SPS.
In this paper, we investigate the factors for the stability and accuracy of the
simulations.
Here, stability means that the solution does not have vibrations in the
simulations, and accuracy means the conservation of constraints in the
simulations.
In general, the accuracy of the simulations would be determined by examining
the numerical solution of the equations.
However, for the nonlinear differential equations, it is often difficult to
investigate the accuracy because of the complexities.
Thus, we adopt the constraints of the system as the criteria of the accuracy in
this paper.

The structure of this paper is as follows.
We review the canonical formulation of the semilinear Klein--Gordon equation in
the de Sitter spacetime in Sec. \ref{sec:CanonicalForm} and the discretized
equation with SPS in Sec. \ref{sec:Discretization}.
In Sec. \ref{sec:Simulation}, we perform some simulations for investigating
their stability and accuracy.
We summarize this paper in Sec. \ref{sec:Summary}.
In this paper, indices such as $(i, j, k, \cdots)$ run from 1 to 3.
We use the Einstein convention of summation of repeated up--down indices.

\section{%
  Canonical formulation of semilinear Klein--Gordon equation in the de Sitter
  spacetime
}%
\label{sec:CanonicalForm}

The semilinear Klein--Gordon equation in the de Sitter spacetime is given by
\begin{eqnarray}
  \partial_t^2\phi
  + 3H\partial_t\phi
  - e^{-2Ht}\delta^{ij}(\partial_i\partial_j \phi)
  + m^2\phi
  + \lambda |\phi|^{p-1}\phi=0,
  \label{eq:KG1}
\end{eqnarray}
where $\phi$ is the field variable, $H$ is the Hubble constant, $\delta^{ij}$
denotes the Kronecker delta, $m$ is the mass, $\lambda$ is a Boolean parameter,
and $p$ is an integer of $2$ or more.
In performing the simulations of Eq. \eqref{eq:KG1}, we often recast first-order
system.
In this paper, we adopt the canonical formulation as the first-order system.
This is because the canonical formulation has the total Hamiltonian, and we can
treat this value as a criterion for investigating the accuracy since the value
is a constraint.

The Hamiltonian density of Eq. \eqref{eq:KG1} is defined as
\begin{eqnarray}
  \mathcal{H}
  &:=&
  \frac{1}{2}e^{-3Ht}\psi^2
  + \frac{1}{2}e^{Ht}\delta^{ij}(\partial_i\phi)(\partial_j\phi)
  + \frac{1}{2}m^2e^{3Ht}\phi^2
  + \frac{\lambda}{p+1}e^{3Ht}|\phi|^{p+1},
\end{eqnarray}
where $\psi$ is the conjugate momentum of $\phi$.
Then, using the canonical equations of $\mathcal{H}$, we obtain the evolution
equations as
\begin{eqnarray}
  \partial_t\phi
  &:=& \frac{\delta\mathcal{H}}{\delta\psi} = e^{-3Ht}\psi,
  \label{eq:phi}
  \\
  \partial_t\psi
  &:=& -\frac{\delta\mathcal{H}}{\delta\phi} =
  e^{Ht}\delta^{ij}(\partial_j\partial_i\phi)
  - m^2e^{3Ht}\phi
  - \lambda e^{3Ht}|\phi|^{p-1}\phi.
  \label{eq:psi}
\end{eqnarray}

The total Hamiltonian $H_C$ is defined as
\begin{eqnarray}
  H_C:=
  \int_{\mathbb{R}^3} \mathcal{H}d^3x,
\end{eqnarray}
and the time derivative of $H_C$ with the evolution Eqs. \eqref{eq:phi} and
\eqref{eq:psi} is
\begin{eqnarray}
  \partial_tH_C
  &=&
  H\int_{\mathbb{R}^3}d^3x\biggl\{
  - \frac{3}{2}e^{-3Ht}\psi^2
  + \frac{1}{2}e^{Ht}\delta^{ij}(\partial_i\phi)(\partial_j\phi)
  + \frac{3}{2}m^2e^{3Ht}\phi^2
  \nonumber\\
  &&
  + \frac{3\lambda}{p+1}e^{3Ht}|\phi|^{p+1}\biggr\}
  + \int_{\mathbb{R}^3}\partial_j\{e^{-3Ht}\delta^{ij}\psi(\partial_i\phi)\}d^3x
  .
  \label{eq:ptHam}
\end{eqnarray}
Note that $H$ is the Hubble constant and $H_C$ is the total Hamiltonian.
If $H=0$ and we set the boundary conditions under which the last term on the
right-hand side of Eq. \eqref{eq:ptHam} is zero on the boundary, then
$\partial_tH_C=0$.
Thus, $H_C$ is treated as a conserved quantity.
On the other hand, in the case of $H\neq 0$, $H_C$ is not a conserved quantity
in general.
In the case of $H\neq 0$, we define the value as
\begin{eqnarray}
  \tilde{H}_C(t)
  &:=& H_C(t) - \int^t_0\partial_sH_C(s)ds.
\end{eqnarray}
$\tilde{H}_C$ identically satisfies $\partial_t\tilde{H}_C=0$.
We call the value $\tilde{H}_C$ as the modified total Hamiltonian hereafter.
In the case of $H\neq 0$, we adopt the value $\tilde{H}_C$ as a criterion for
the accuracy of the simulations.
To investigate the accuracy of the simulations, we monitor $H_C$ in a flat
spacetime such as $H=0$, and $\tilde{H}_C$ in a nonflat spacetime such as
$H=10^{-3}$.
If the changes in $H_C$ in the flat spacetime or $\tilde{H}_C$ in the nonflat
spacetime against the initial values are sufficiently small during the
evolution, we determine that the simulations are successful.
That is, the smaller the change in $H_C$ or $\tilde{H}_C$ in the evolution, the
more accurate the numerical calculations.

\section{%
  Discretizations of semilinear Klein--Gordon equation in the de
  Sitter spacetime
}%
\label{sec:Discretization}
The main factor for the numerical errors occurs during the processes of the
discretizations of the equations.
In this section, we review the discretized equations of the semilinear
Klein--Gordon equation in the de Sitter spacetime.

The discretized Hamiltonian density is defined as
\begin{eqnarray}
  \mathcal{H}^{(\ell)}_{(\boldsymbol{k})}
  &:=&
  \frac{1}{2}e^{-3Ht^{(\ell)}}(\psi^{(\ell)}_{(\boldsymbol{k})})^2
  + \frac{1}{2}e^{Ht^{(\ell)}}\delta^{ij}
  (\hat{\delta}^{\langle1\rangle}_i\phi^{(\ell)}_{(\boldsymbol{k})})
  (\hat{\delta}^{\langle1\rangle}_j\phi^{(\ell)}_{(\boldsymbol{k})})
  \nonumber\\
  &&
  + \frac{1}{2}m^2e^{3Ht^{(\ell)}}(\phi^{(\ell)}_{(\boldsymbol{k})})^2
  + \frac{\lambda}{p+1}e^{3Ht^{(\ell)}}|\phi^{(\ell)}_{(\boldsymbol{k})}|^{p+1}.
  \label{eq:DiscreteHam}
\end{eqnarray}
By using SPS, we can rewrite the discretized Eqs. \eqref{eq:phi} and
\eqref{eq:psi} as
\begin{eqnarray}
  \frac{\phi^{(\ell+1)}_{(\boldsymbol{k})}
    - \phi^{(\ell)}_{(\boldsymbol{k})}}{\Delta t}
  &=&\frac{1}{4}(e^{-3Ht^{(\ell+1)}} + e^{-3Ht^{(\ell)}})
  (\psi^{(\ell+1)}_{(\boldsymbol{k})} + \psi^{(\ell)}_{(\boldsymbol{k})}),
  \label{eq:DiscretePhi}
  \\
  \frac{\psi^{(\ell+1)}_{(\boldsymbol{k})}
    - \psi^{(\ell)}_{(\boldsymbol{k})}}{\Delta t}
  &=&
  \frac{1}{4}(e^{Ht^{(\ell+1)}}+e^{Ht^{(\ell)}})
  \delta^{ij}\widehat{\delta}^{\langle1\rangle}_i
  \widehat{\delta}^{\langle1\rangle}_j(\phi^{(\ell+1)}_{(\boldsymbol{k})}
  + \phi^{(\ell)}_{(\boldsymbol{k})})
  \nonumber\\
  &&
  - \frac{m^2}{4}(e^{3Ht^{(\ell+1)}} + e^{3Ht^{(\ell)}})
  (\phi^{(\ell+1)}_{(\boldsymbol{k})} + \phi^{(\ell)}_{(\boldsymbol{k})})
  \nonumber\\
  &&
  - \frac{\lambda }{2(p+1)}(e^{3Ht^{(\ell+1)}} +e^{3Ht^{(\ell)}})
  \frac{|\phi^{(\ell+1)}_{(\boldsymbol{k})}|^{p+1}
    - |\phi^{(\ell)}_{(\boldsymbol{k})}
    |^{p+1}}{\phi^{(\ell+1)}_{(\boldsymbol{k})}
    - \phi^{(\ell)}_{(\boldsymbol{k})}},
  \label{eq:DiscretePsi}
\end{eqnarray}
respectively.
The upper index ${}^{(\ell)}$ in parentheses is the time index, and the lower
index ${}_{(\boldsymbol{k})}$ in parentheses is the spatial grid index, where
$\boldsymbol{k}=(k_1,k_2,k_3)$ and $k_1$, $k_2$, and $k_3$ are
$x$, $y$, and $z$ indices, respectively.
$\widehat{\delta}^{\langle1\rangle}_i$ is the discrete operator defined as
\begin{eqnarray}
  \widehat{\delta}^{\langle1\rangle}_i u^{(\ell)}_{(\boldsymbol{k})}
  := \left\{
  \begin{array}{ll}
    \dfrac{u^{(\ell)}_{(k_1+1,k_2,k_3)} - u^{(\ell)}_{(k_1-1,k_2,k_3)}}{
      2\Delta x},
    & (i=1)\\
    \dfrac{u^{(\ell)}_{(k_1,k_2+1,k_3)} - u^{(\ell)}_{(k_1,k_2-1,k_3)}}{
      2\Delta y},
    & (i=2)\\
    \dfrac{u^{(\ell)}_{(k_1,k_2,k_3+1)} - u^{(\ell)}_{(k_1,k_2,k_3-1)}}{
      2\Delta z}.
    & (i=3)\\
  \end{array}
  \right.
\end{eqnarray}
There are two features in Eq. \eqref{eq:DiscretePsi}.
First, the second-order difference is expressed as
$\widehat{\delta}^{\langle1\rangle}_i\widehat{\delta}^{\langle1\rangle}_j$.
In general, the discrete operator of the second-order difference is usually
defined as
\begin{eqnarray}
  \widehat{\delta}^{\langle2\rangle}_{ij}u^{(\ell)}_{(\boldsymbol{k})}
  := \left\{
  \begin{array}{ll}
    \dfrac{u^{(\ell)}_{(k_1+1,k_2,k_3)} - 2u^{(\ell)}_{(\boldsymbol{k})}
      + u^{(\ell)}_{(k_1-1,k_2,k_3)}}{(\Delta x)^2},
    & (i=j=1)\\
    \dfrac{u^{(\ell)}_{(k_1,k_2+1,k_3)} - 2u^{(\ell)}_{(\boldsymbol{k})}
      + u^{(\ell)}_{(k_1,k_2-1,k_3)}}{(\Delta y)^2},
    & (i=j=2)\\
    \dfrac{u^{(\ell)}_{(k_1,k_2,k_3+1)} - 2u^{(\ell)}_{(\boldsymbol{k})}
      + u^{(\ell)}_{(k_1,k_2,k_3-1)}}{(\Delta z)^2},
    & (i=j=3)\\
    \widehat{\delta}^{\langle1\rangle}_i\widehat{\delta}^{\langle1\rangle}_j
    u^{(\ell)}_{(\boldsymbol{k})}.
    & (i\neq j)\\
  \end{array}
  \right.
\end{eqnarray}
In the case of $i=j$, $\widehat{\delta}^{\langle2\rangle}_{ij}u^{(\ell)}_{
  (\boldsymbol{k})}
\neq \widehat{\delta}^{\langle1\rangle}_{i}
\widehat{\delta}^{\langle1\rangle}_{j}u^{(\ell)}_{(\boldsymbol{k})}$.
Second, the expression of the nonlinear term, which is the last term on the
right-hand side in Eq. \eqref{eq:DiscretePsi}, is not usual.
In general, the discretized expression expected from Eq. \eqref{eq:psi} is
$-\lambda e^{3Ht^{(\ell)}}|\phi^{(\ell)}_{(\boldsymbol{k})}|^{p-1}
\phi^{(\ell)}_{(\boldsymbol{k})}$.
These differences in the simulations are shown in Sec. \ref{sec:Simulation}.

The discretized total Hamiltonian $H_C^{(\ell)}$ is defined as
\begin{eqnarray}
  H^{(\ell)}_{C}
  &:=&
  \sum_{\substack{1\le k_1\le n_1\\
      1\le k_2\le n_2\\
      1\le k_3\le n_3}}
  \mathcal{H}^{(\ell)}_{(\boldsymbol{k})}
  \Delta x\Delta y\Delta z,
  \label{eq:disTotalHam}
\end{eqnarray}
where $n_1$, $n_2$, and $n_3$ are the grid numbers for $x$, $y$, and $z$,
respectively.
The difference quotient for $H_C^{(\ell)}$ using Eqs. \eqref{eq:DiscretePhi}
and \eqref{eq:DiscretePsi} is calculated as
\begin{eqnarray}
  &&\frac{H^{(\ell+1)}_C- H^{(\ell)}_C}{\Delta t}
  \nonumber\\
  &&=
  H\sum_{\substack{1\le k_1\le n_1\\
      1\le k_2\le n_2\\
      1\le k_3\le n_3}}\biggl[
    -\frac{3}{4}\{e^{-3Ht^{(\ell+1)}}(\psi^{(\ell+1)}_{(\boldsymbol{k})})^2
    + e^{-3Ht^{(\ell)}}(\psi^{(\ell)}_{(\boldsymbol{k})})^2\}
  \nonumber\\
  &&\quad
  + \frac{1}{4}\delta^{ij}\{
  e^{Ht^{(\ell+1)}}(\hat{\delta}^{\langle1\rangle}_i\phi^{(\ell+1)}_{(\boldsymbol{k})})
  (\hat{\delta}^{\langle1\rangle}_j\phi^{(\ell+1)}_{(\boldsymbol{k})})
  + e^{Ht^{(\ell)}}(\hat{\delta}^{\langle1\rangle}_i\phi^{(\ell)}_{(\boldsymbol{k})})
  (\hat{\delta}^{\langle1\rangle}_j\phi^{(\ell)}_{(\boldsymbol{k})})
  \}
  \nonumber\\
  &&\quad
  + \frac{3}{4}m^2\{e^{3Ht^{(\ell+1)}}(\phi^{(\ell+1)}_{(\boldsymbol{k})})^2
  + e^{3Ht^{(\ell)}}(\phi^{(\ell)}_{(\boldsymbol{k})})^2\}
  \nonumber\\
  &&\quad
  + \frac{3\lambda}{2(p+1)}
  (e^{3Ht^{(\ell+1)}}|\phi^{(\ell+1)}_{(\boldsymbol{k})}|^{p+1}
  + e^{3Ht^{(\ell)}}|\phi^{(\ell)}_{(\boldsymbol{k})}|^{p+1})
  \biggr]
  \nonumber\\
  &&\quad
  +[\text{Boundary\,\, Terms}]
  + O(\Delta t),
  \label{eq:discHCprop1}
\end{eqnarray}
where we use the relation such that
\begin{eqnarray}
  e^{at^{(\ell+1)}}
  = e^{at^{(\ell)}}
  + ae^{at^{(\ell)}}\Delta t + O((\Delta t)^2).
  \quad (\forall a\in\mathbb{R})
\end{eqnarray}
The boundary terms in Eq. \eqref{eq:discHCprop1} are eliminated under the
periodic boundary condition.
In addition, if $H=0$, then $H^{(\ell+1)}_C$ is consistent with $H^{(0)}_C$
in the order of $\Delta t$.
Then we define the discretized modified total Hamiltonian
$\tilde{H}^{(\ell)}_C$ as
\begin{eqnarray}
  \tilde{H}^{(\ell)}_{C}
  &:=&
  H^{(\ell)}_C
  - H\sum_{0\leq m\leq \ell-1}
  \sum_{\substack{1\le k_1\le n_1\\
      1\le k_2\le n_2\\
      1\le k_3\le n_3}}\biggl[
    -\frac{3}{4}\{e^{-3Ht^{(m+1)}}(\psi^{(m+1)}_{(\boldsymbol{k})})^2
    + e^{-3Ht^{(m)}}(\psi^{(m)}_{(\boldsymbol{k})})^2\}
    \nonumber\\
    &&
    + \frac{1}{4}\delta^{ij}\{
    e^{Ht^{(m+1)}}(\hat{\delta}^{\langle1\rangle}_i\phi^{(m+1)}_{(\boldsymbol{k})})
    (\hat{\delta}^{\langle1\rangle}_j\phi^{(m+1)}_{(\boldsymbol{k})})
    + e^{Ht^{(m)}}(\hat{\delta}^{\langle1\rangle}_i\phi^{(m)}_{(\boldsymbol{k})})
    (\hat{\delta}^{\langle1\rangle}_j\phi^{(m)}_{(\boldsymbol{k})})
    \}
    \nonumber\\
    &&
    + \frac{3}{4}m^2\{e^{3Ht^{(m+1)}}(\phi^{(m+1)}_{(\boldsymbol{k})})^2
    + e^{3Ht^{(m)}}(\phi^{(m)}_{(\boldsymbol{k})})^2\}
    \nonumber\\
    &&
    + \frac{3\lambda}{2(p+1)}
    (e^{3Ht^{(m+1)}}|\phi^{(m+1)}_{(\boldsymbol{k})}|^{p+1}
    + e^{3Ht^{(m)}}|\phi^{(m)}_{(\boldsymbol{k})}|^{p+1})
    \biggr]
  \Delta t\Delta x\Delta y\Delta z.
  \label{eq:tildeHam1}
\end{eqnarray}
We adopt this value as a criterion of the accuracy of the simulations in the
case of $H\neq 0$.
\section{Numerical simulations}
\label{sec:Simulation}

In this section, we perform some simulations with SPS to investigate their
stability and accuracy.
We perform simulations with three formulations of the discretized semilinear
Klein--Gordon equation in the de Sitter spacetime.
The first formulation is that for Eqs. \eqref{eq:DiscretePhi},
\eqref{eq:DiscretePsi}, and \eqref{eq:disTotalHam}.
We call this formulation Form I.
As shown in Eq. \eqref{eq:discHCprop1}, Form I is SPS.
The details are shown in \cite{Tsuchiya-Nakamura-2019-JCAM}.
The second formulation is that for Eqs. \eqref{eq:DiscretePhi},
\eqref{eq:disTotalHam}, and the following Eq. \eqref{eq:DiscretePsi2}.
\begin{eqnarray}
  \frac{\psi^{(\ell+1)}_{(\boldsymbol{k})} - \psi^{(\ell)}_{(\boldsymbol{k})}}{
    \Delta t}
  &=&
  \frac{1}{4}(e^{Ht^{(\ell+1)}}+e^{Ht^{(\ell)}})
    \delta^{ij}\widehat{\delta}^{\langle2\rangle}_{ij}
    (\phi^{(\ell+1)}_{(\boldsymbol{k})} + \phi^{(\ell)}_{(\boldsymbol{k})})
    \nonumber\\
  &&
    - \frac{m^2}{4}(e^{3Ht^{(\ell+1)}} + e^{3Ht^{(\ell)}})
    (\phi^{(\ell+1)}_{(\boldsymbol{k})} + \phi^{(\ell)}_{(\boldsymbol{k})})
    \nonumber\\
  &&
    - \frac{\lambda }{2(p+1)}(e^{3Ht^{(\ell+1)}} +e^{3Ht^{(\ell)}})
  \frac{|\phi^{(\ell+1)}_{(\boldsymbol{k})}|^{p+1}
    - |\phi^{(\ell)}_{(\boldsymbol{k})}|^{p+1}}{
    \phi^{(\ell+1)}_{(\boldsymbol{k})}
    - \phi^{(\ell)}_{(\boldsymbol{k})}}
  \label{eq:DiscretePsi2}
\end{eqnarray}
We call this formulation Form II.
The difference between Eqs. \eqref{eq:DiscretePsi} and \eqref{eq:DiscretePsi2}
is the second-order difference term.
The third formulation is that for Eqs. \eqref{eq:DiscretePhi},
\eqref{eq:disTotalHam}, and the following Eq. \eqref{eq:DiscretePsi3}.
\begin{eqnarray}
  \frac{\psi^{(\ell+1)}_{(\boldsymbol{k})} - \psi^{(\ell)}_{(\boldsymbol{k})}}{
    \Delta t}
  &&=
  \frac{1}{4}(e^{Ht^{(\ell+1)}}+e^{Ht^{(\ell)}})
  \delta^{ij}\widehat{\delta}^{\langle1\rangle}_{i}
  \widehat{\delta}^{\langle1\rangle}_{j}
  (\phi^{(\ell+1)}_{(\boldsymbol{k})} + \phi^{(\ell)}_{(\boldsymbol{k})})
  \nonumber\\
  &&
  - \frac{m^2}{4}(e^{3Ht^{(\ell+1)}} + e^{3Ht^{(\ell)}})
  (\phi^{(\ell+1)}_{(\boldsymbol{k})} + \phi^{(\ell)}_{(\boldsymbol{k})})
  \nonumber\\
  &&
  - \frac{\lambda }{8}(e^{3Ht^{(\ell+1)}} +e^{3Ht^{(\ell)}})
  |\phi^{(\ell+1)}_{(\boldsymbol{k})}+\phi^{(\ell)}_{(\boldsymbol{k})}|^{p-1}
  (\phi^{(\ell+1)}_{(\boldsymbol{k})} + \phi^{(\ell)}_{(\boldsymbol{k})})
  \label{eq:DiscretePsi3}
\end{eqnarray}
We call this formulation Form III.
The difference between Eqs. \eqref{eq:DiscretePsi} and \eqref{eq:DiscretePsi3}
is the expression of the discretized nonlinear term, which is the last term on
the right-hand side of each of these equations.

The simulation settings are as follows.
\begin{itemize}
\item Initial conditions: $\phi_0=A\cos(2\pi x)$, $\psi_0=2\pi A\sin(2\pi x)$,
  and $A=4$
\item Numerical domains: $0\leq x\leq 1$, $0\leq t\leq 1000$
\item Boundary condition: periodic
\item Grids: $\Delta x=1/200$ and $\Delta t=1/1000$
\item Mass: $m=1$
\item Boolean parameter of the nonlinear term: $\lambda=1$
\item Number of exponents in the nonlinear term: $p=2,3,4,5$, and $6$
\item Hubble constant: $H=0$ and  $10^{-3}$
\end{itemize}
Forms I, II, and III are expressed in three dimensions.
On the other hand, the initial conditions are one-dimensional.
Even if the spatial dimension of the initial conditions is one-dimensional, the
differences exist in the second-order difference term and the discretized
nonlinear term.
Thus, the numerical simulations are expected to show differences in the
one-dimensional initial conditions.

\subsection{%
  Flat spacetime
}%
\label{subsec:flat}

\begin{figure}[t]
  \centering
  \includegraphics[width=0.32\hsize]{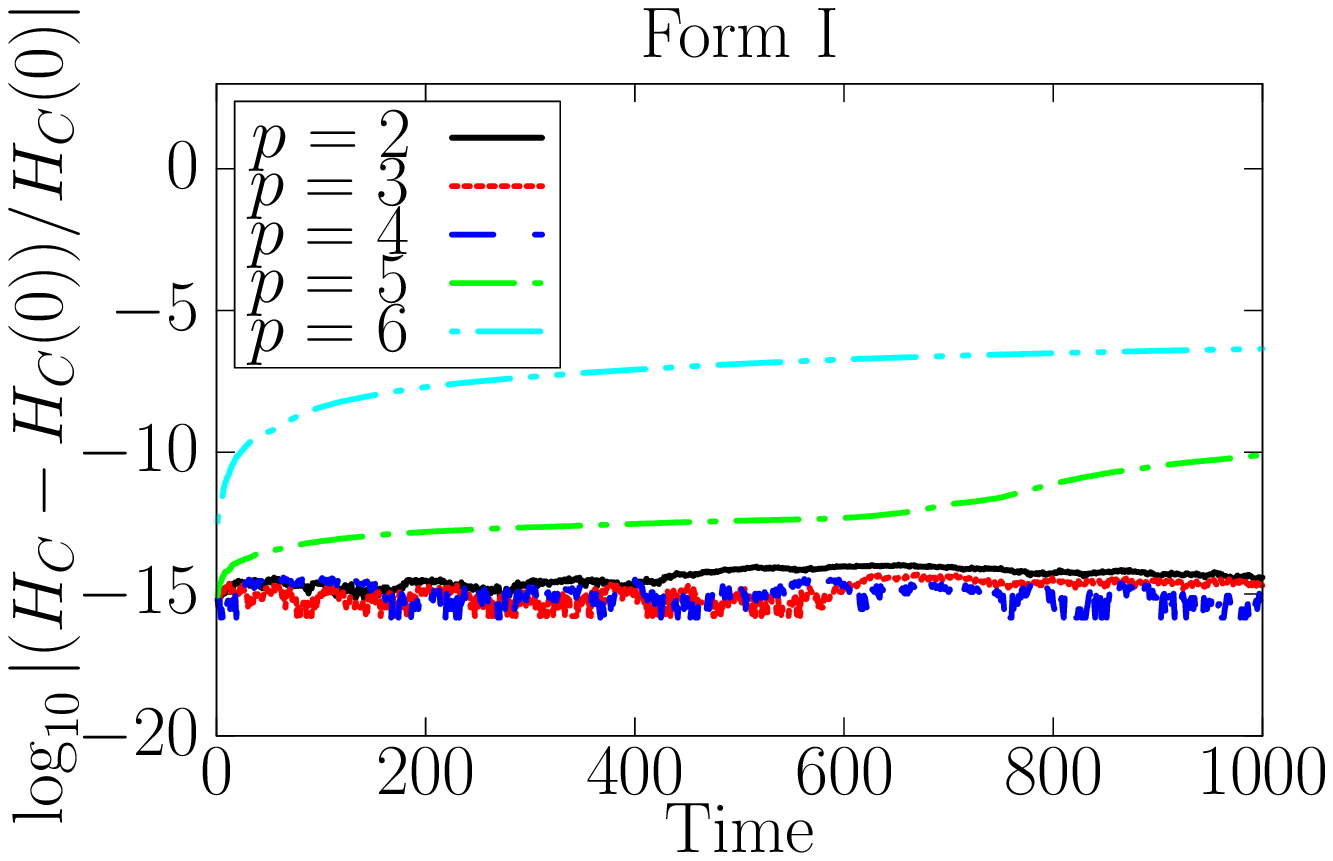}
  \includegraphics[width=0.32\hsize]{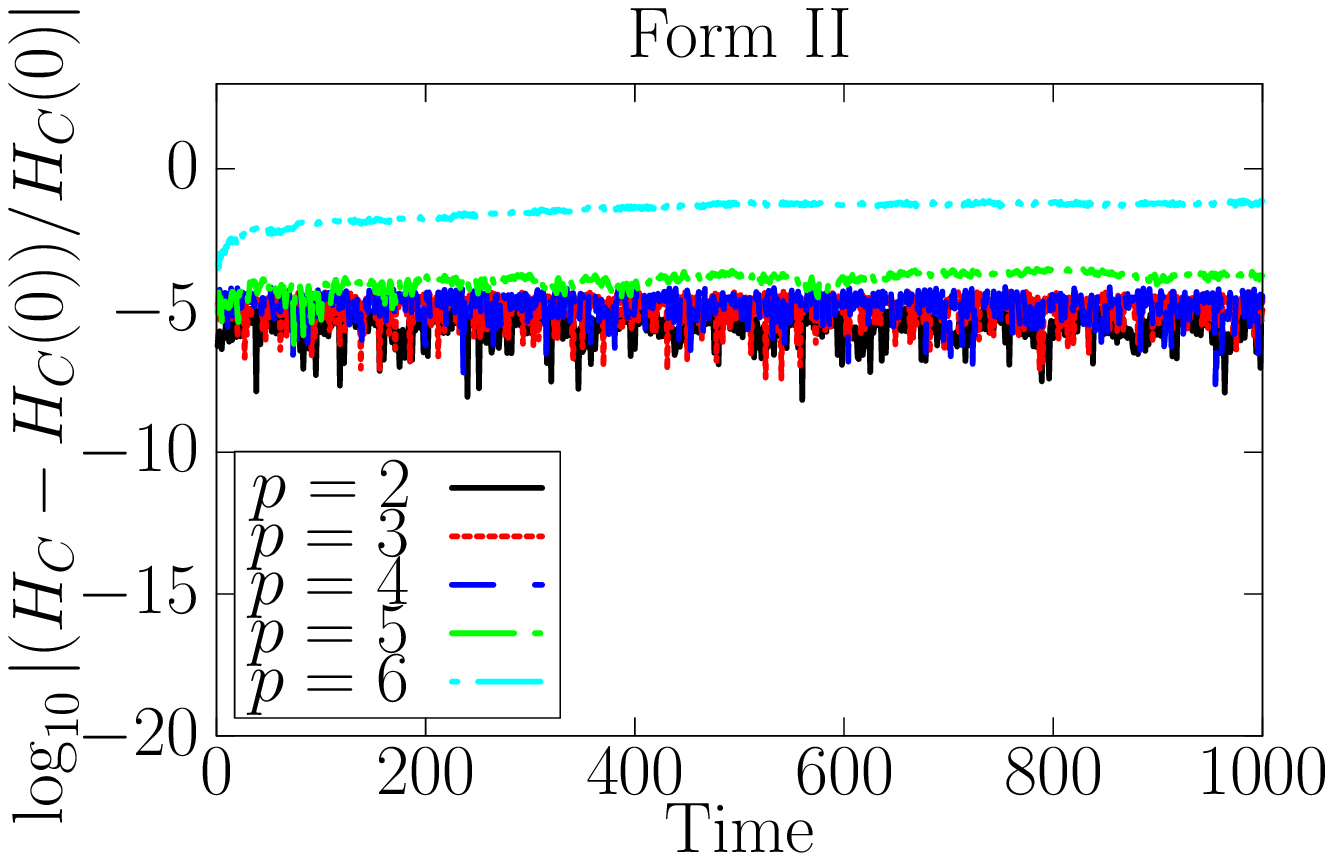}
  \includegraphics[width=0.32\hsize]{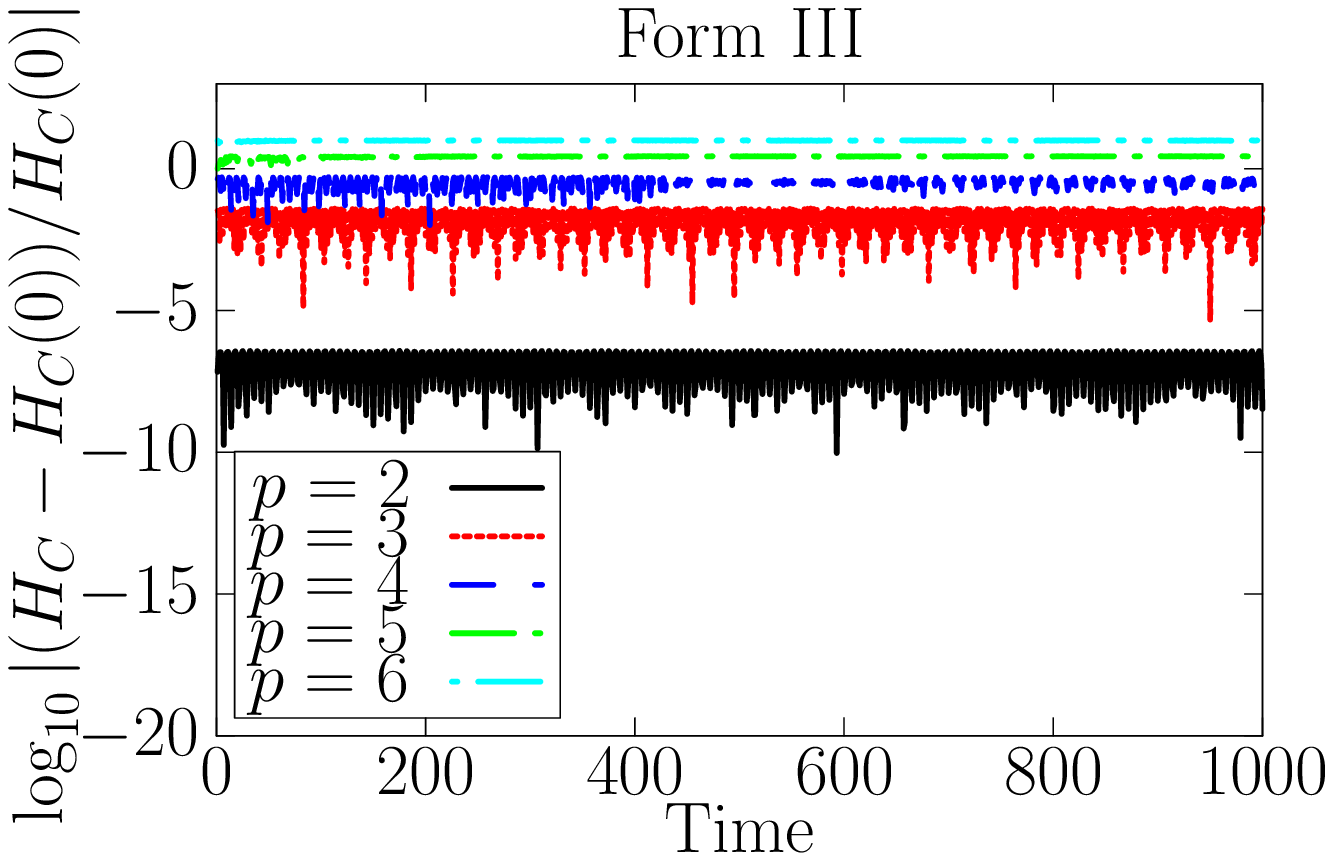}
  \caption{%
    Relative errors of the total Hamiltonian $H_C$ against the initial value
    $H_C(0)$ for each $p$ in the case of $H=0$.
    The horizontal axis is time, and the vertical axis is
    $\log_{10}|(H_C-H_C(0))/H_C(0)|$.
    The left panel is drawn with Form I, the center panel with Form II, and the
    right panel with Form III.
  }
  \label{fig:totalH_flat}
\end{figure}

We perform some simulations of the three formulations in the flat spacetime,
which is in the case of $H=0$.
In Fig. \ref{fig:totalH_flat}, we show the relative errors of the total
Hamiltonian $H_C$ against the initial values $H_C(0)$ for each value of the
exponent $p$ in the nonlinear term.
The left panel is drawn with Form I, the center panel with Form II, and the
right panel with Form III.
The values of $|(H_C-H_C(0))/H_C(0)|$ indicate the numerical errors because
$H_C$ is a constraint.
In the right panel, we see that the value of $p=2$ with Form III is smaller
than those of the other exponents in the panel.
This result indicates that the numerical errors caused by the nonlinear term are
small in the case of $p=2$.
We see that the values of the center and right panels are larger than that of
the left panel for each $p$.
Thus, the simulations with Form I are more accurate than those with the other
forms.

\begin{figure}[t]
  \centering
  \includegraphics[width=0.32\hsize]{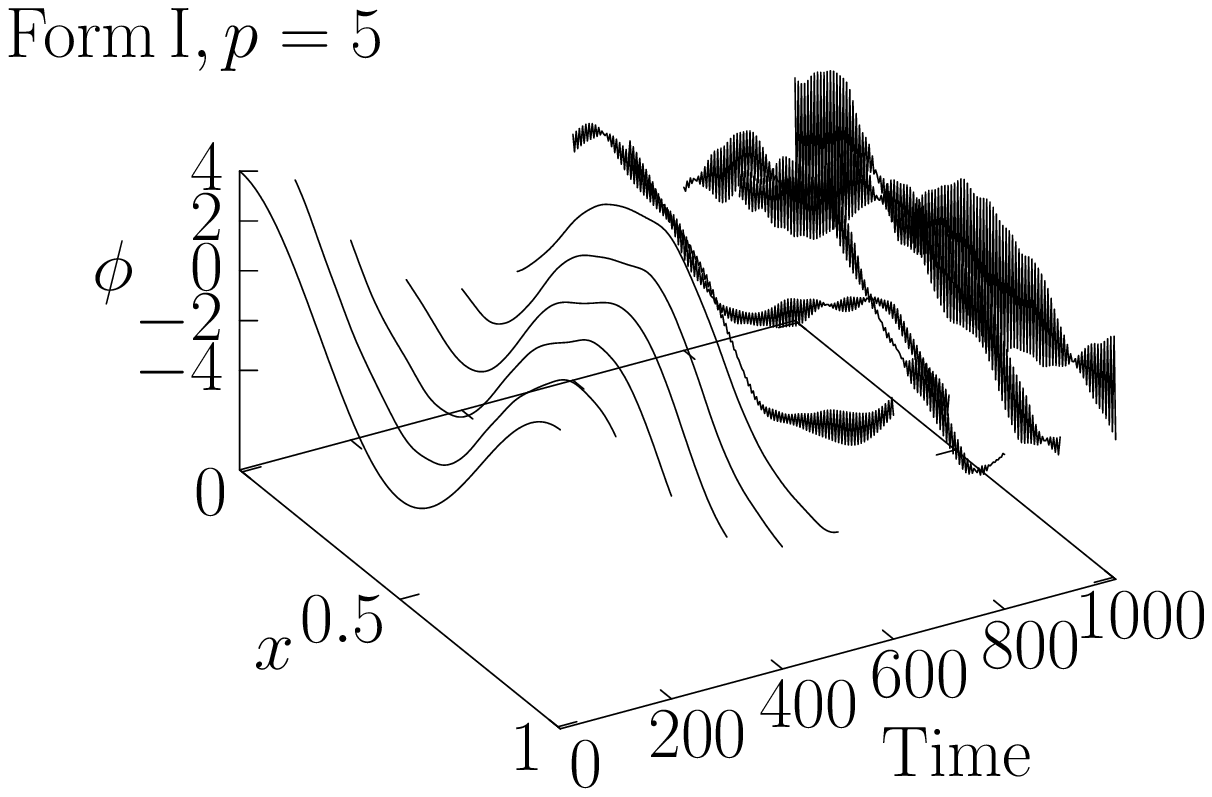}
  \includegraphics[width=0.32\hsize]{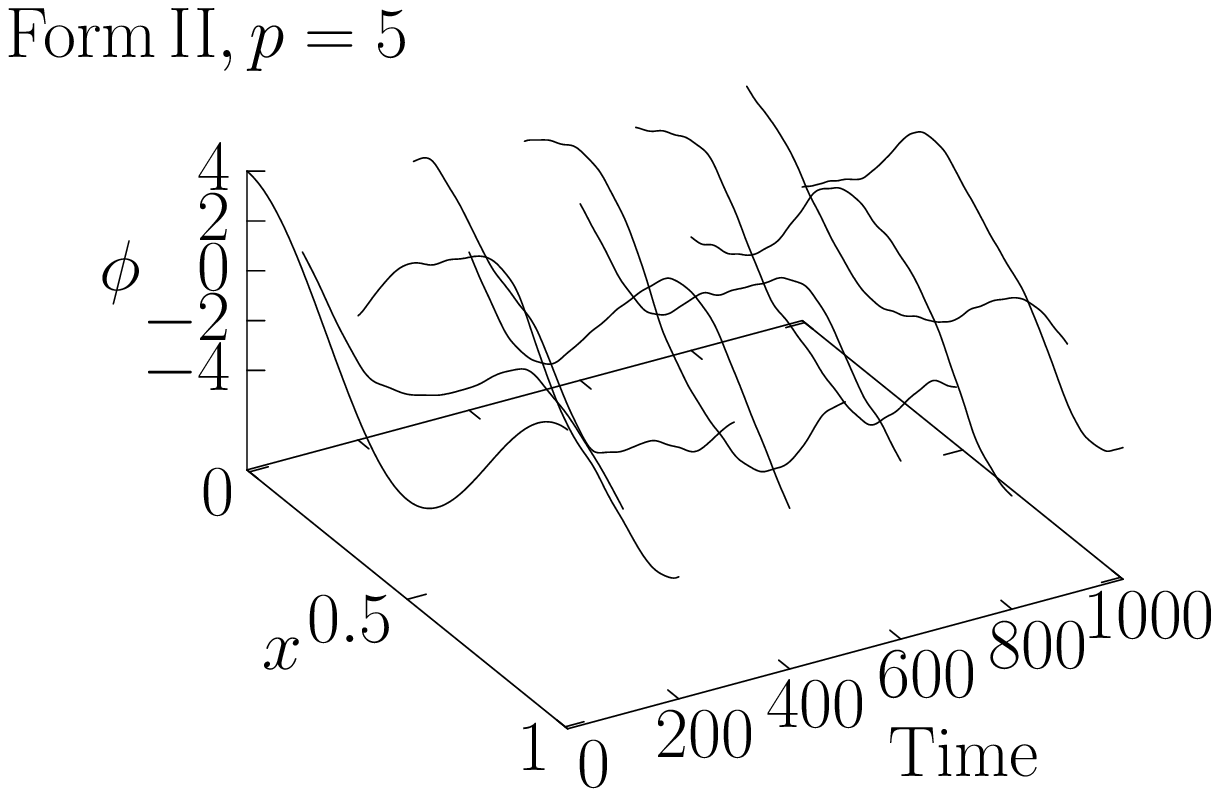}
  \includegraphics[width=0.32\hsize]{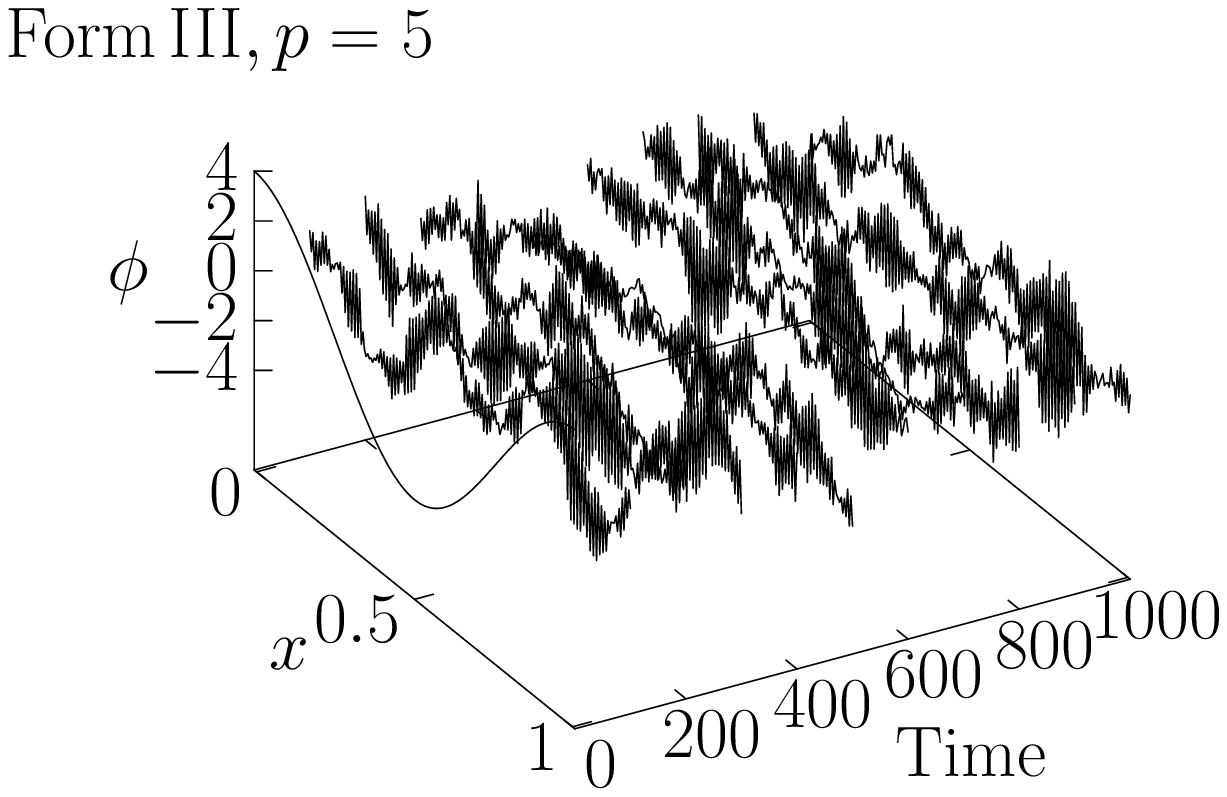}\\
  \includegraphics[width=0.32\hsize]{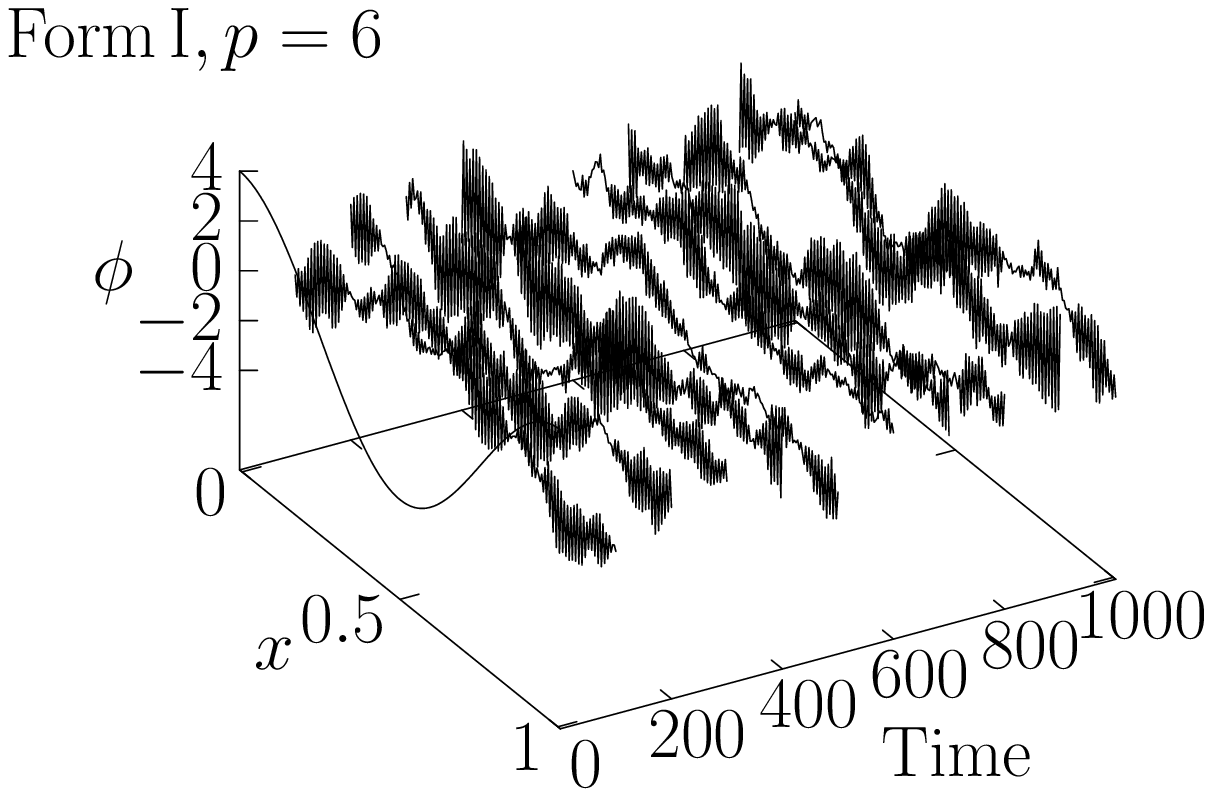}
  \includegraphics[width=0.32\hsize]{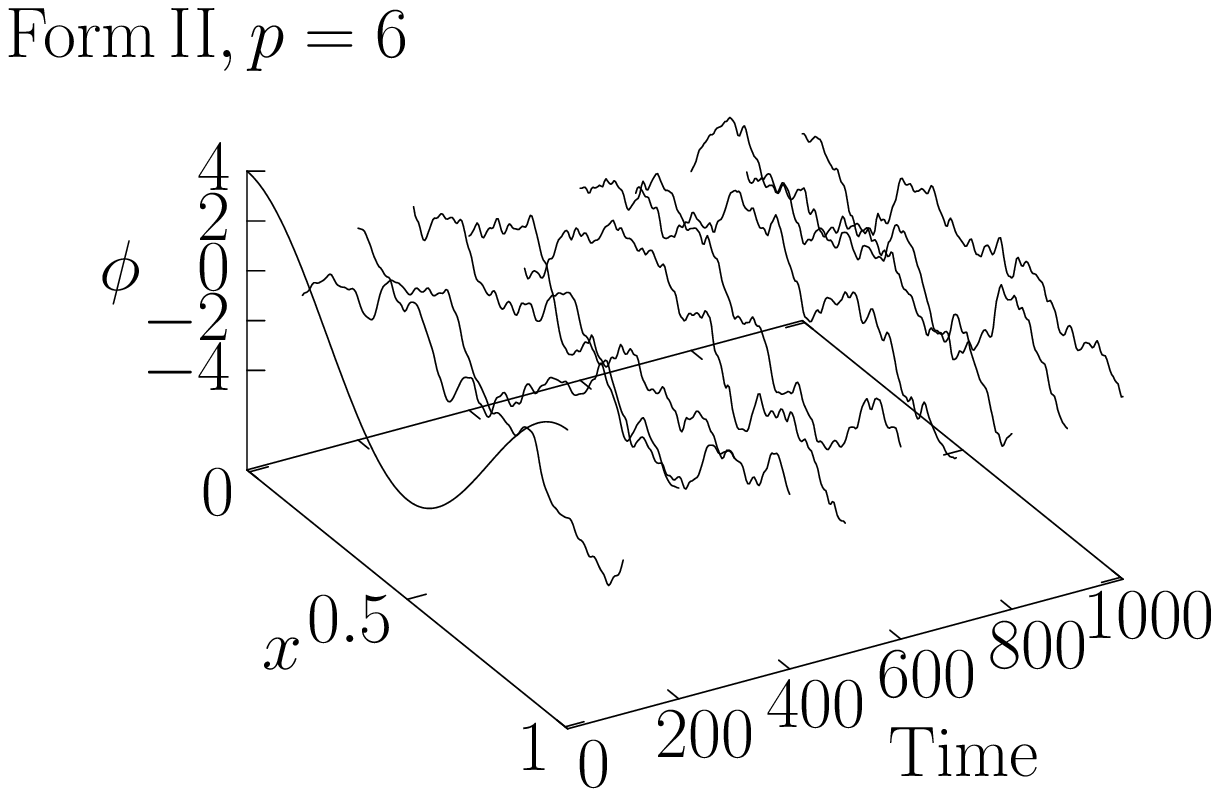}
  \includegraphics[width=0.32\hsize]{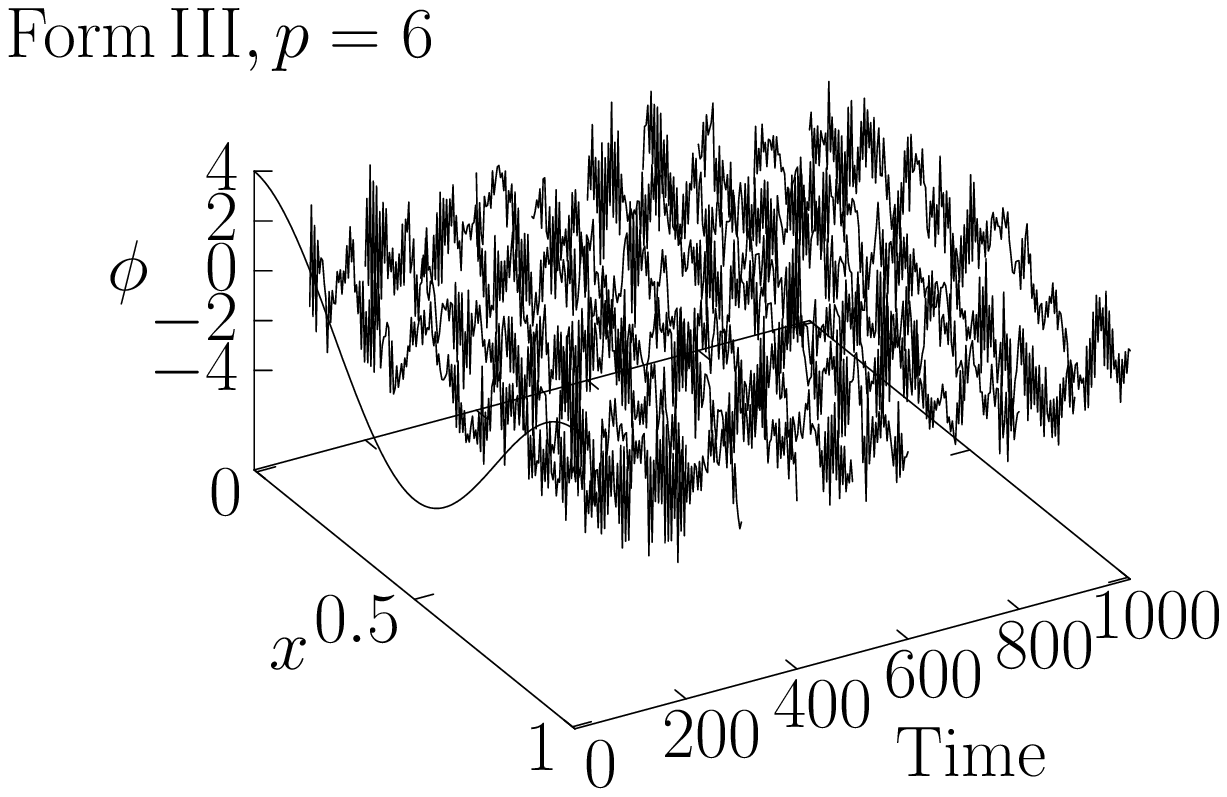}
  \caption{%
    $\phi$ with $p=5$ and $6$.
    The left panels are drawn with Form I, the center panels with Form II,
    and the right panels with Form III.
    The top panels are drawn for $p=5$ and the bottom panels for $p=6$.
    The vibrations occur at $t\geq 700$ in the top-left panel, $t\geq 100$
    in the bottom-left panel, and $t\geq 100$ in the right panels.
  }
  \label{fig:phi_flat}
\end{figure}

Then we show $\phi$ with $p=5$ and $6$ in Fig. \ref{fig:phi_flat} to investigate
the stability of the simulations.
The left panels are drawn with Form I, the center panels with Form II, and the
right panels with Form III.
The top panels are drawn with the exponent $p=5$ and the bottom panels with
$p=6$.
We see that the simulations of the top-left panel at $t\geq 700$, the
bottom-left panel at $t\geq 100$, and the right panels at $t\geq 100$ are
unstable because of the generated vibrations.
On the other hand, the simulations shown in the center panels are stable until
$t=1000$.
Thus, we determine that the simulations with Form II are more stable than those
with the other formulations.

\subsection{%
  Curved spacetime
}%
\label{subsec:curved}

Here, we perform some simulations with the same settings as in Sec.
\ref{subsec:flat} except for the Hubble constant.
This time, we set the Hubble constant $H=10^{-3}$.

\begin{figure}[t]
  \centering
  \includegraphics[width=0.32\hsize]{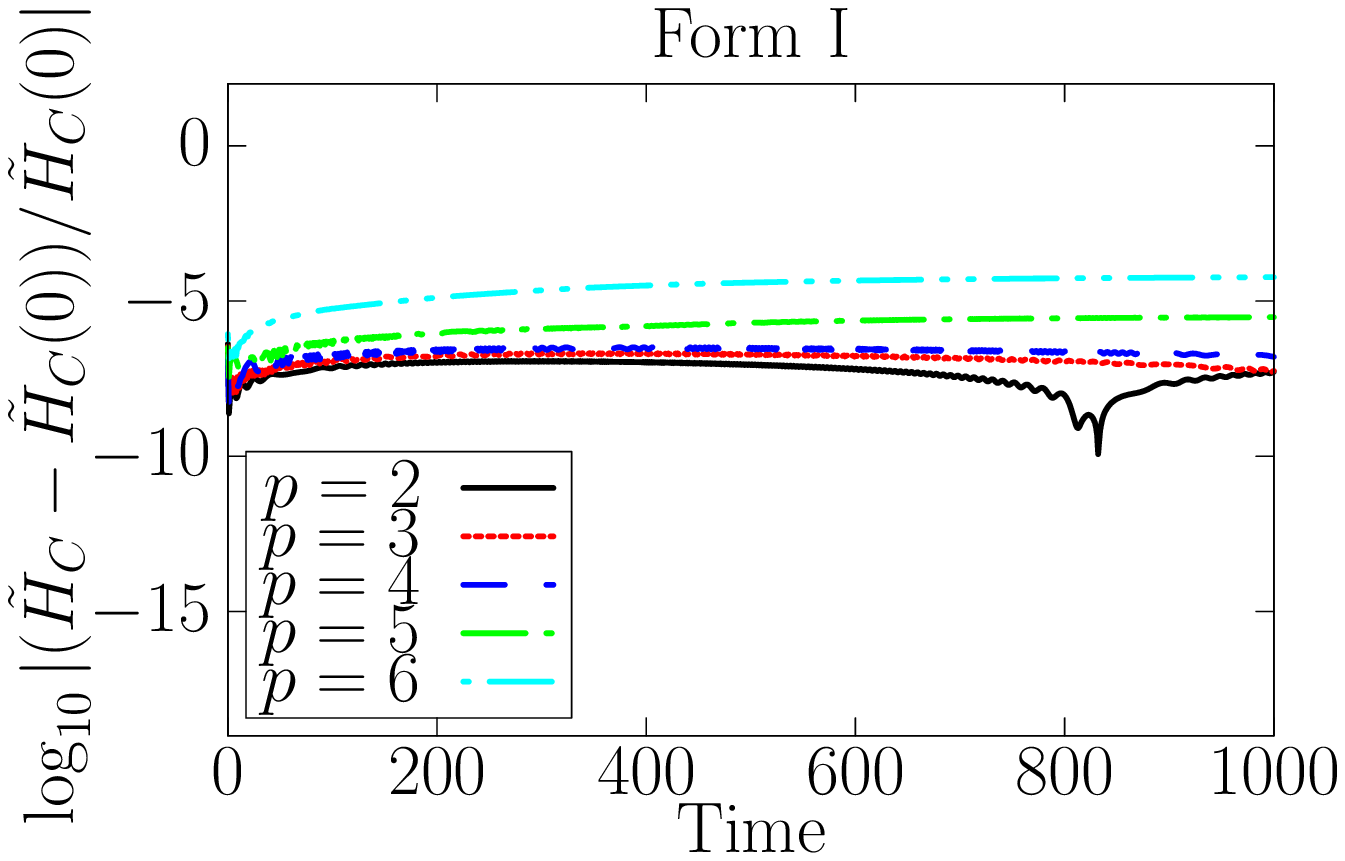}
  \includegraphics[width=0.32\hsize]{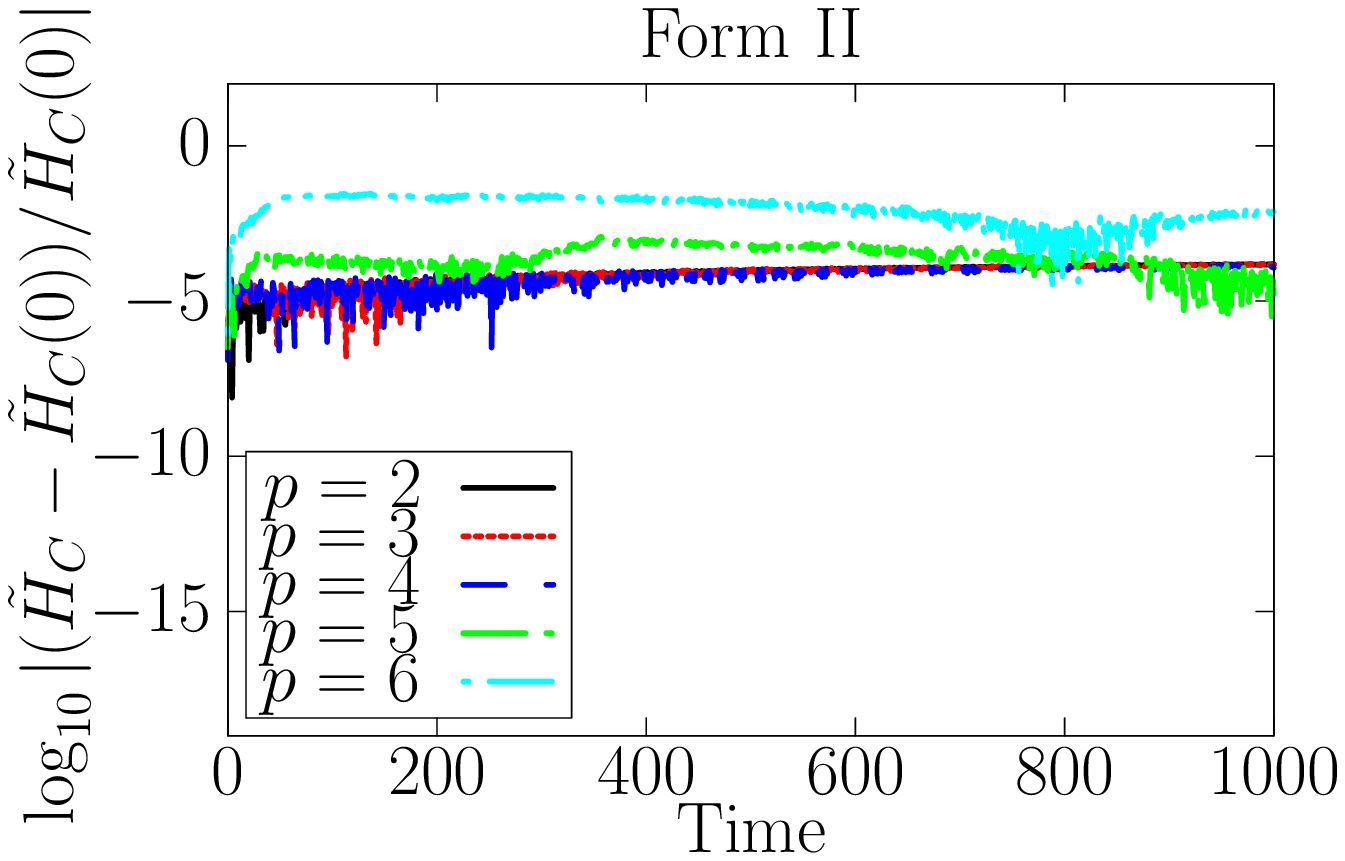}
  \includegraphics[width=0.32\hsize]{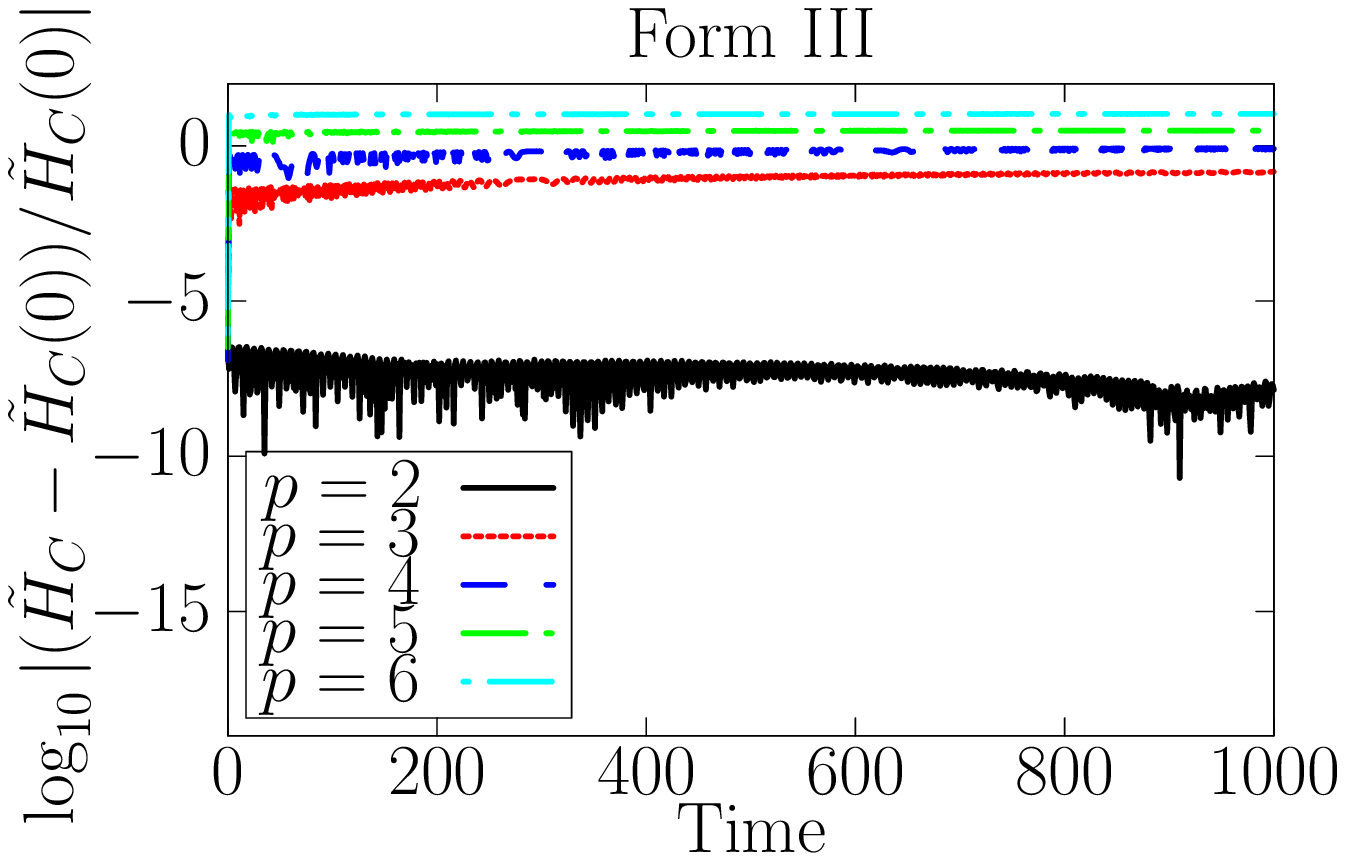}
  \caption{%
    Relative errors of the modified total Hamiltonian $\tilde{H}_C$ against
    the initial value $\tilde{H}_C(0)$ for each $p$ in the case of $H=10^{-3}$.
    The horizontal axis is time, and the vertical axis is
    $\log_{10}|(\tilde{H}_C-\tilde{H}_C(0))/\tilde{H}_C(0)|$.
    The left panel is drawn with Form I, the center panel with Form II, and the
    right panel with Form III.
  }
  \label{fig:modHC}
\end{figure}

We show the relative errors of the modified total Hamiltonian $\tilde{H}_C$
against the initial value $\tilde{H}_C(0)$ in Fig. \ref{fig:modHC}.
Note that $\tilde{H}_C$ is calculated approximately using Eq.
\eqref{eq:tildeHam1} via the numerical solutions in time evolution.
The left panel is drawn with Form I, the center panel with Form II, and the
right panel with Form III.
We see that the value of $p=2$ with Form III is smaller than those in the other
cases in the right panel.
This tendency is consistent with the case of $H=0$.

\begin{figure}[t]
  \centering
  \includegraphics[width=0.32\hsize]{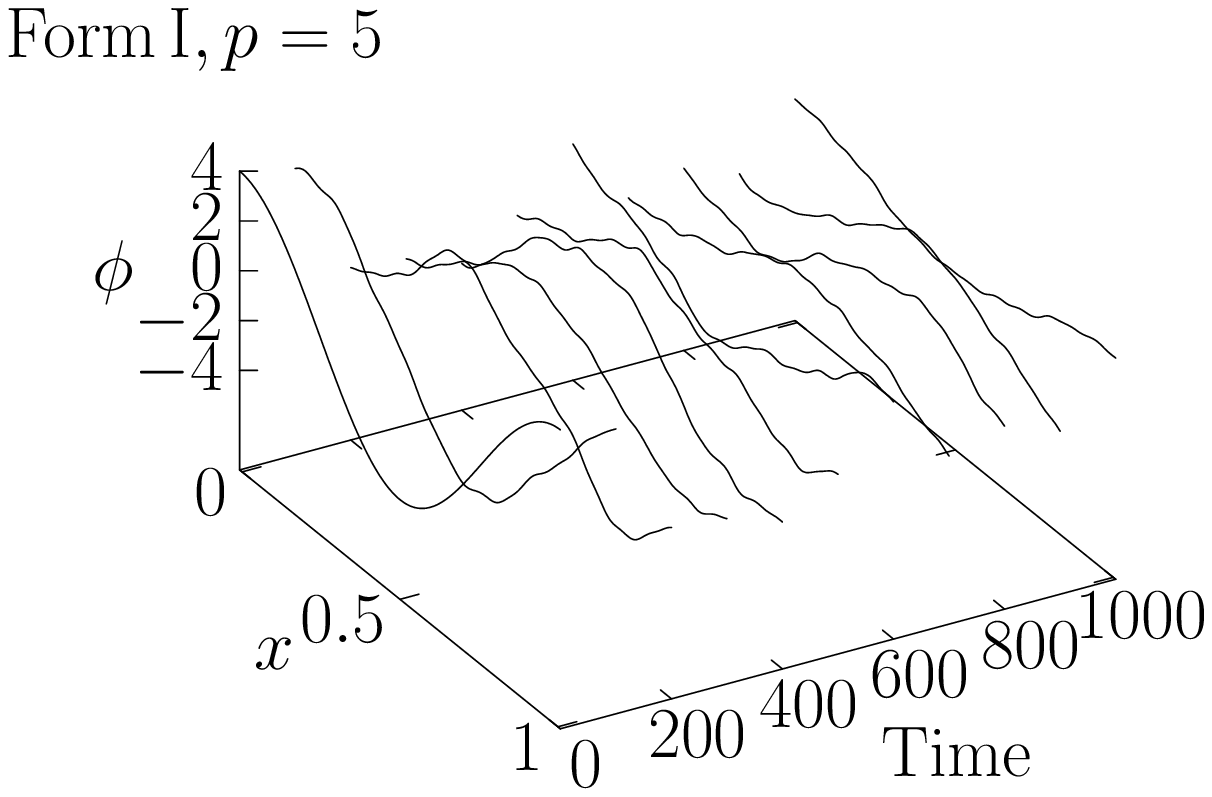}
  \includegraphics[width=0.32\hsize]{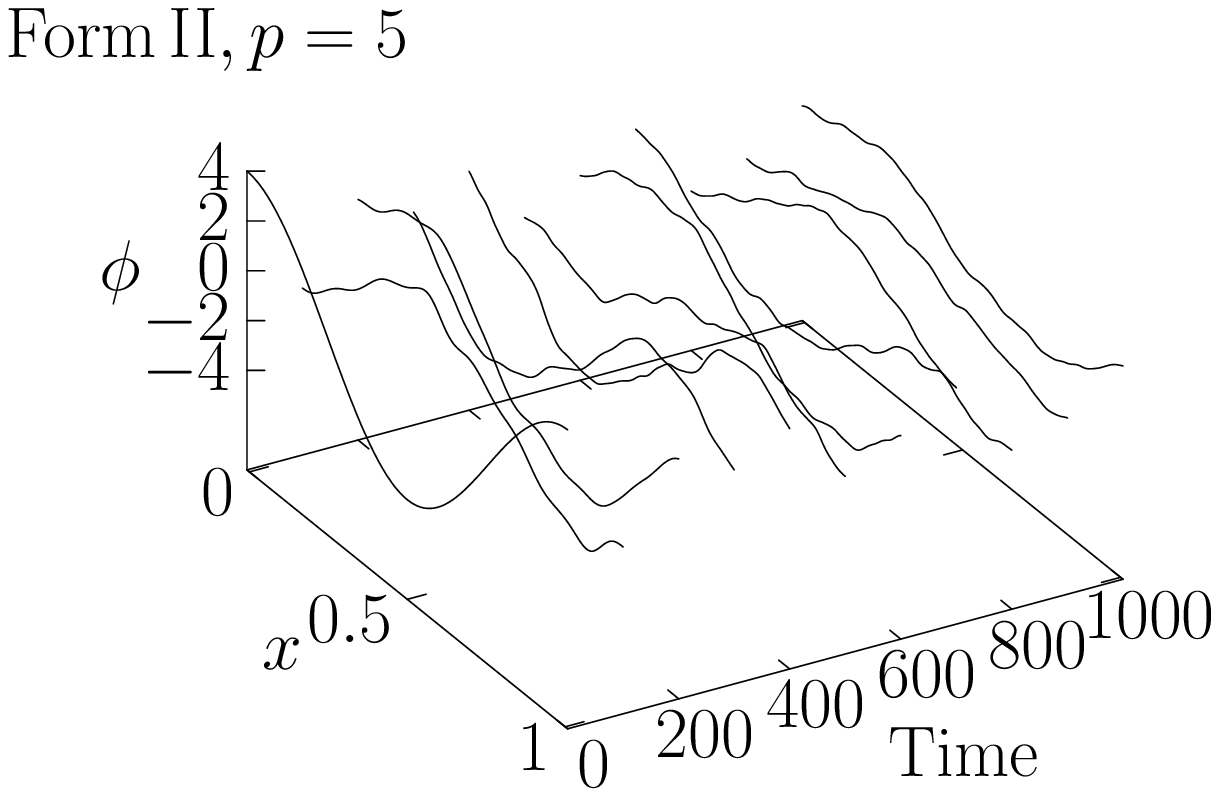}
  \includegraphics[width=0.32\hsize]{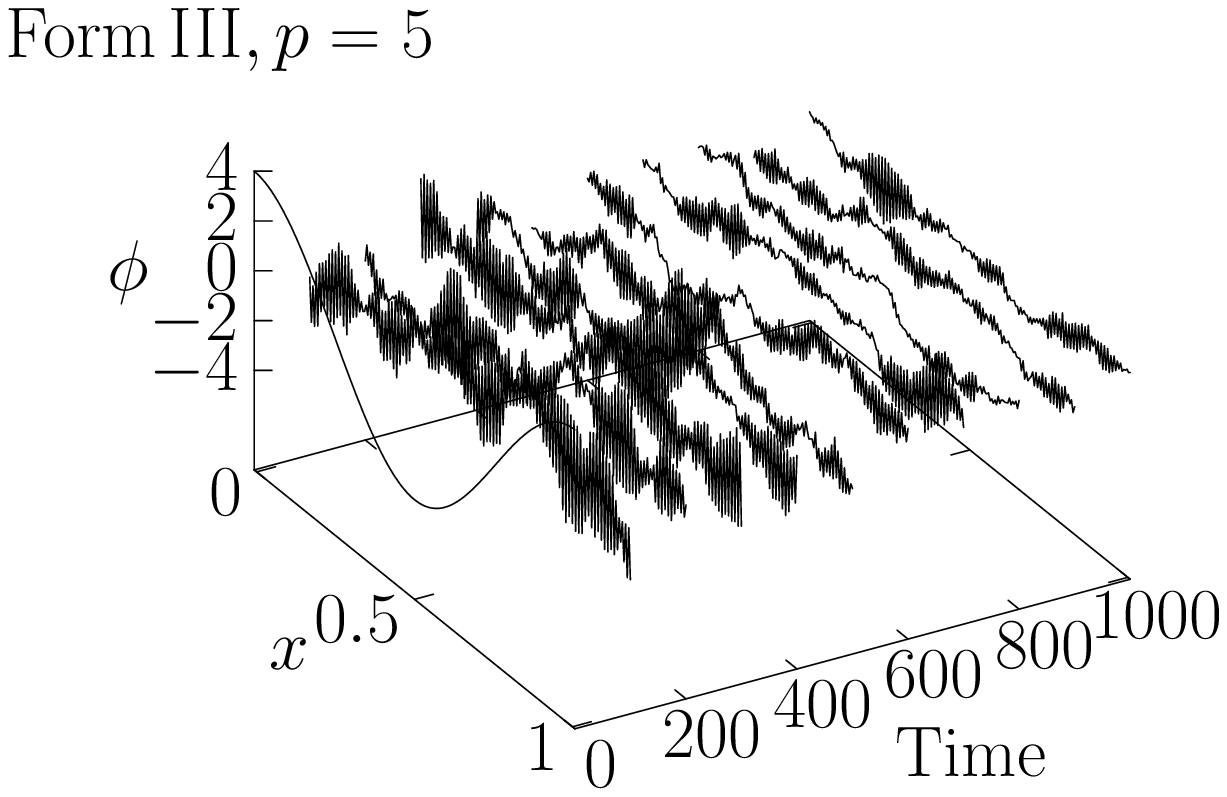}\\
  \includegraphics[width=0.32\hsize]{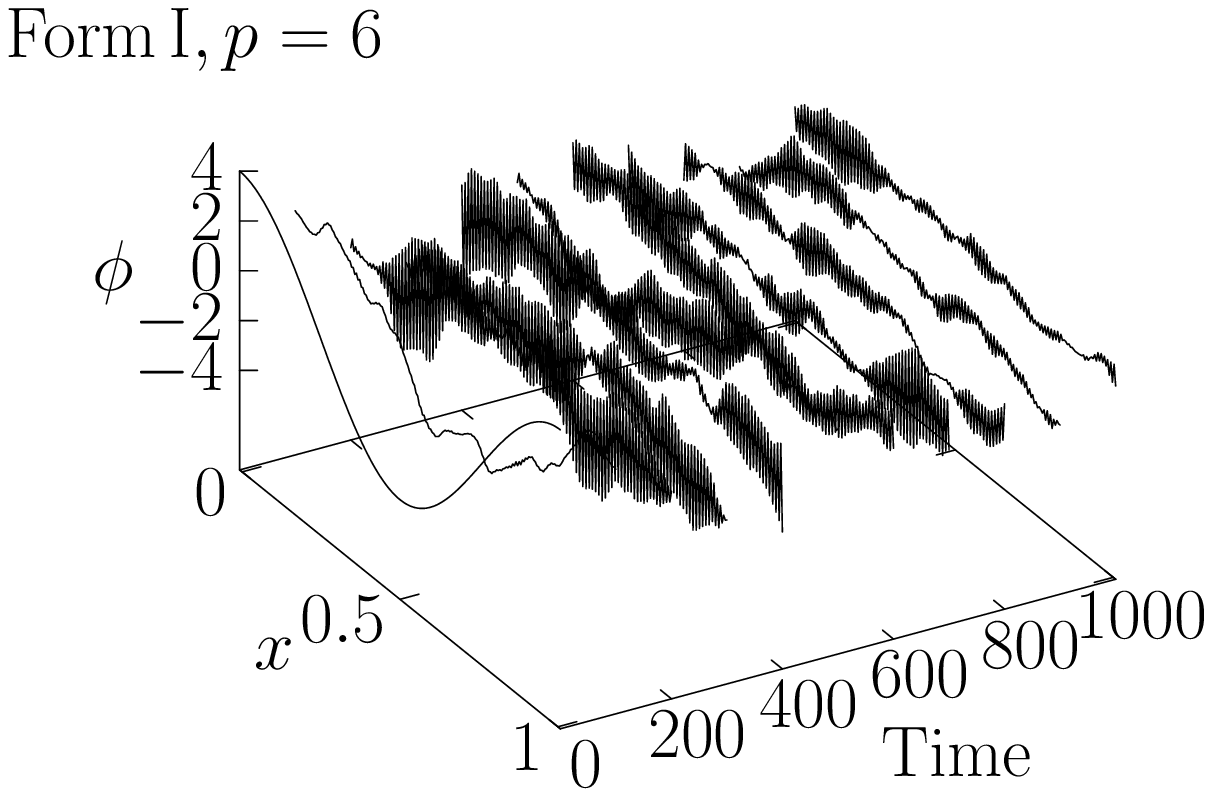}
  \includegraphics[width=0.32\hsize]{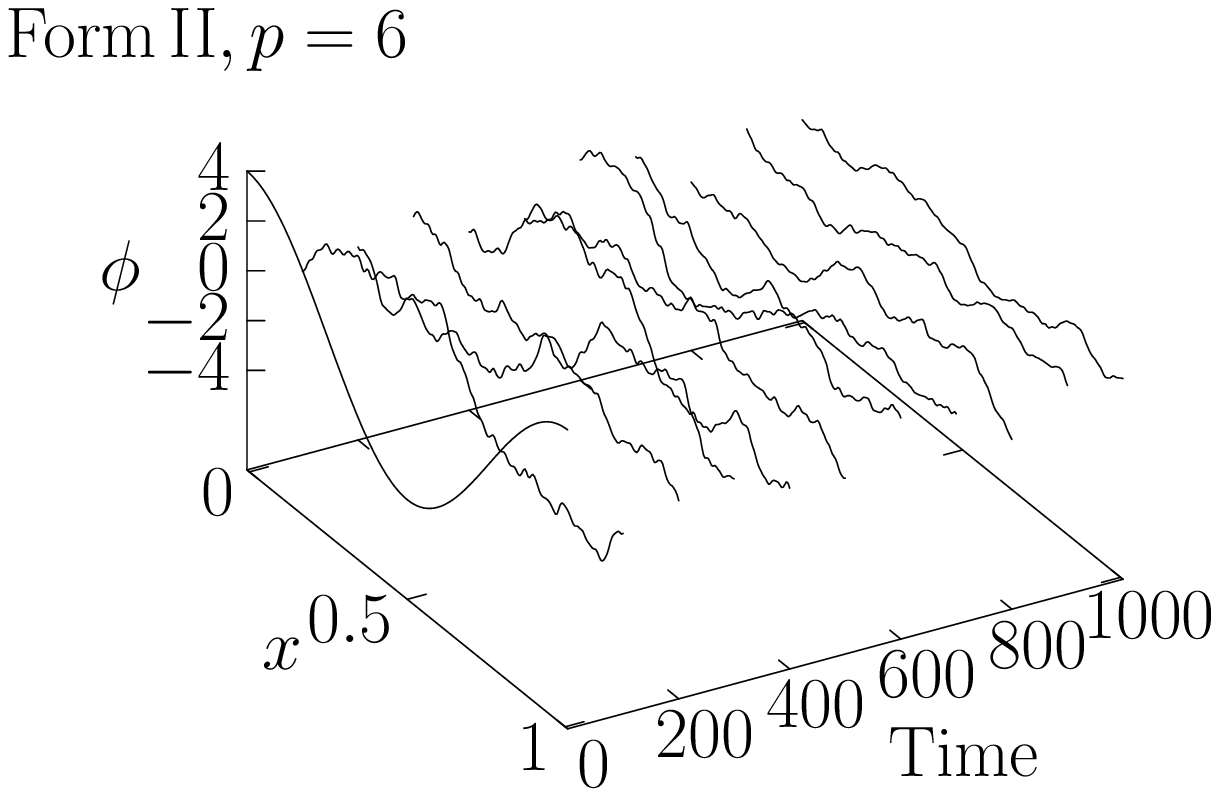}
  \includegraphics[width=0.32\hsize]{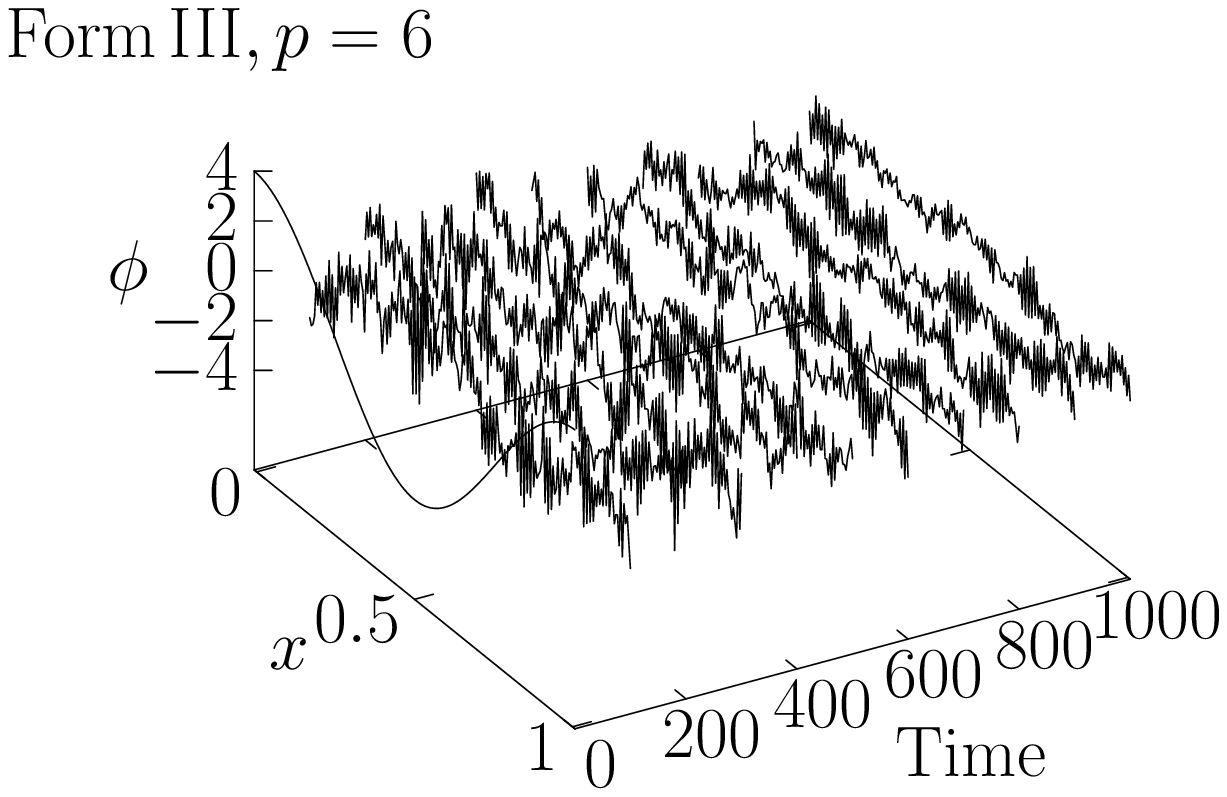}
  \caption{
    The same as in Fig. \ref{fig:phi_flat} except for the value of the Hubble
    constant, which is $10^{-3}$.
  }
  \label{fig:phi_curved}
\end{figure}

Figure \ref{fig:phi_curved} is the same as Fig. \ref{fig:phi_flat} except for
the value of the Hubble constant.
In the comparison between Figs. \ref{fig:phi_flat} and \ref{fig:phi_curved}, no
vibrations appear in the top-left panel in Fig. \ref{fig:phi_curved}
and also in the bottom-left panel in Fig. \ref{fig:phi_curved} until $t=100$.
The other patterns of behavior are almost the same.
These results indicate that the vibrations of the waveform of the solutions
decrease in comparison with the case of $H=0$.
That is, the positive Hubble constant makes the simulation stable.
This is also noted in \cite{Tsuchiya-Nakamura-2019-JCAM}.

\section{Summary}
\label{sec:Summary}

We investigated the factors affecting the stability and accuracy of simulations
of the semilinear Klein--Gordon equation in the de Sitter spacetime using SPS.
We reviewed the canonical formulation of the equation and that of the
discretized equation with SPS.
To investigate the terms affecting the stability and accuracy in the
discretized equations, we compared some simulations using three discretized
formulations.
The first formulation consists of the discretized equations with SPS, which is
called Form I.
This formulation was reported in \cite{Tsuchiya-Nakamura-2019-JCAM}.
The second formulation consists of the discretized equations with SPS, in which
the second-order difference was replaced with a standard discretized
second-order difference, which is called Form II.
The third formulation consists of the discretized equations with SPS, in which
the nonlinear term was replaced with a standard discretized term, which is
called Form III.
We monitored the total Hamiltonian or the modified one to see the accuracy of
the simulations.
As a result, we found that the stability and accuracy of the simulations using
Form III are worse than those with Form I.
This result indicates that the discretizations of the nonlinear term affect on
the stability and accuracy of the simulations.
In addition, the accuracy of the simulations with Form I is
better than those with the other forms.
On the other hand, the stability of the simulations with Form II is higher than
those with the other forms.
Moreover, we confirmed that the simulations with positive values of the Hubble
constant are more stable than those in the flat spacetime.

The numerical stability of the simulations using Form I is lower than those
using Form II.
However, there are degrees of freedom in the selection of the discretized terms
for Form I.
Therefore, it seems that the formulation that enables stable and accurate
numerical simulation can be constructed, which we will report in the near
future.

\section*{Acknowledgments}
The authors thank the anonymous referees for their many helpful comments
that improved the paper.
T.T. and M.N. were partially supported by JSPS KAKENHI Grant Number 21K03354.
T.T. was partially supported by JSPS KAKENHI Grant Number 20K03740
and Grant for Basic Science Research Projects from The Sumitomo Foundation.
M.N. was partially supported by JSPS KAKENHI Grant Number 16H03940.


\end{document}